\documentclass[11pt,reqno]{amsart}
\usepackage{amsmath,amssymb,amsthm,upref,graphicx,mathrsfs}
\usepackage{color,xcolor}
\usepackage[
  colorlinks=true,
  linkcolor=blue,
  citecolor=blue,
  urlcolor=blue]{hyperref}
\usepackage{xfrac}

\newtheorem{thm}{Theorem}[section]
\newtheorem{cor}[thm]{Corollary}
\newtheorem{lem}[thm]{Lemma}
\newtheorem{prop}[thm]{Proposition}

\newtheorem{conj}[thm]{Conjecture}
\numberwithin{equation}{section}

\begin{document}

\leftline{ \scriptsize}

\vspace{1.3 cm}
\title
{Idempotent factorization on some matrices over quadratic integer rings}
\author{ Peeraphat Gatephan$^{\ast}$ and Kijti Rodtes$^{\ast\ast}$ }
\thanks{{\scriptsize
			$^{*}$ Department of Mathematics, Faculty of Science, Naresuan University, Phitsanulok 65000, Thailand, E-mail addresses: peeraphatg64@nu.ac.th (Peeraphat Gatephan).\\
			$^{**}$ Department of Mathematics and Research center for Academic Excellence in Mathematics,  Faculty of Science, Naresuan University, Phitsanulok 65000, Thailand, E-mail addresses:  kijtir@nu.ac.th (Kijti Rodtes).\\
		\\}}
\hskip -0.4 true cm

\maketitle


\begin{abstract} In 2020, Cossu and Zanardo raised a conjecture on the idempotent factorization on singular matrices in the form $\begin{pmatrix}
		p&z\\ \bar{z}&\sfrac{\lVert z\rVert}{p}
	\end{pmatrix},$ where $p$ is a prime integer which is irreducible but not prime element in the ring of integers $\mathbb{Z}[\sqrt{D}]$ and $z\in\mathbb{Z}[\sqrt{D}]$  such that  $\langle p,z\rangle$ is a non-principal ideal.  In this paper, we provide some classes of matrices that affirm the conjecture and some classes of matrices that oppose the conjecture.  We further show that there are matrices in the above form that can not be written as a product of  two idempotent matrices.
\end{abstract}

\vskip 0.2 true cm


\pagestyle{myheadings}
\markboth{\rightline {\scriptsize Peeraphat Gatephan and Kijti Rodtes}}
{\leftline{\scriptsize }}
\bigskip
\bigskip


\vskip 0.4 true cm

\section{Introduction}

In 1967, J.A. Erdos \cite{e67} proved that every singular square matrix over fields admits idempotent factorizations. This result leads to a characterization of integral domains $R$ satisfying the property $\operatorname{ID}_n$: every $n\times n$ singular matrix over $R$ can be written as a product of idempotent matrices over $R$. The idempotent factorizations problem motivated Salce and Zanardo \cite{sz14} to conjecture that an integral domain satisfying $\operatorname{ID}_2$ must be a B\'{e}zout domain (every finitely generated ideal is principal). 

 In 2020, Cossu and Zanardo \cite{cz20} showed that every singular matrix having at least one row or one column whose elements generate a principal ideal must admit $\operatorname{ID}_2$ (see, Corollary 4.1. in \cite{cz20} for details). However, the complete classification of $2\times2$ singular matrices over ring of integers of quadratic number field $\mathbb{Z}[\sqrt{D}]$ ($D$ is a square free integer) is still an open problem. In addition, they investigated the factorization on some matrices in the form $$A(p,z):=\begin{pmatrix}
 p&z\\ \bar{z}&\sfrac{\lVert z\rVert}{p}
 \end{pmatrix}$$   where $p$ is a prime integer which is irreducible but not prime element in $\mathbb{Z}[\sqrt{D}]$ and $z\in\mathbb{Z}[\sqrt{D}]$  such that the ideal generated by $p,z$,  $\langle p,z\rangle$, is a non-principal ideal.  They also raised the following conjecture in the same paper.
\begin{conj}\label{con1}\cite{cz20} Let $D$ be a square free integer, $p$ be a prime integer which is irreducible but not prime in $\mathbb{Z}[\sqrt{D}]$ and $z\in\mathbb{Z}[\sqrt{D}]$ be such that $\langle p,z\rangle$ is a non-principal ideal. Then
	$$A(p,z)=\begin{pmatrix}
		a&b\\ c&1-a
	\end{pmatrix}\begin{pmatrix}
		\bar{a}&\bar{c}\\ \bar{b}& 1-\bar{a}
	\end{pmatrix},$$ for some $a,b,c\in\mathbb{Z}[\sqrt{D}]$ with $a(1-a)=bc$.
\end{conj}

In this paper, we apply row operations on a system of Diophantine equations over $\mathbb{Q}[\sqrt{D}]$ to get a necessary and sufficient condition for the factorization in term of quadratic Diophantine equations.  According to this condition, we can conclude  in Theorem \ref{thmnz=p^2} that if $\lVert z\rVert=-p^2$ then the matrix $A(p,z)$ always affirms the Cossu and Zanardo's conjecture.   Moreover, by using this condition together with the Florida Transform and modulo trick, there are explicit classes of matrices $A(p,z)$ that oppose the conjecture in  Corollary \ref{corflorida} and Corollary \ref{coroppose}.

\section{Background and basic results}

Let $R$ be an integral domains. For any $a_1,a_{2},\ldots,a_k\in R$, we denote $\langle a_1,a_{2},\ldots,a_k\rangle$ the ideal generated by $a_1,a_{2},\ldots,a_k$. A principal ideal is an ideal in a ring $R$ that is generated by a single element of $R$. We use $R^{\times }$ to represent the group of unit elements in $R$ and $\mathbb{M}_n(R)$ to represent the set of all $n\times n$ matrices over $R$. For any $A\in\mathbb{M}_n(R)$, $A$ is idempotent if and only if $A^2=A$.  By direct computation, it is known that a singular nonzero matrix $A=(A_{ij})\in\mathbb{M}_2(R)$ is idempotent if and only if $A_{22}=1-A_{11}$. For any $A\in\mathbb{M}_n(R)$, we say that $A\in\operatorname{ID}_n(R)$ if $A$ can be expressed as a product of idempotent matrices in $\mathbb{M}_n(R)$ . 

Throughout this paper, we always denote $D$ a square free integer, $p$ a prime number  which is  irreducible but not prime  in $\mathbb{Z}[\sqrt{D}]$.  For a given $D$ and $p$ as above, a characterization of  $z\in\mathbb{Z}[\sqrt{D}]$ for which $\langle p, z\rangle$ is non-principal can be determined via the set $I_p(D)$: \textit{the set of all non unit $z\in\mathbb{Z}[\sqrt{D}]$ for which $z\notin\langle p\rangle$ but there exists $m\notin\langle p\rangle$  such that $zm\in\langle p\rangle.$}  Precisely, 
\begin{prop}\label{propnz}
For any $z\in \mathbb{Z}[\sqrt{D}]$,  we have
	\begin{center}
		$z\in I_p(D)$ if and only if $\langle p,z\rangle$ is a non-principal ideal.
	\end{center}
\end{prop}
\begin{proof}
	Let $z\in I_p(D)$ and assume for a contradiction that $\langle p,z\rangle=\langle c\rangle$ for some $c\in\mathbb{Z}[\sqrt{D}]$.   Since $z\in I_p(D)$, there exists $m\notin\langle p\rangle $ such that $zm=pl$ for some $m,l\in\mathbb{Z}[\sqrt{D}]$.   Since $p\in\langle p,z\rangle=\langle c\rangle$, there exists $a\in\mathbb{Z}[\sqrt{D}]$ such that $p=ca$.  By using the fact that $p$ is an irreducible element, we have $c$ or $a$ must be unit in $\mathbb{Z}[\sqrt{D}]$.  If $c$ is a unit, then $\langle c\rangle=\mathbb{Z}[\sqrt{D}]$. So, $1\in \langle c\rangle=\langle p,z\rangle$; namely, there exists $x,y\in\mathbb{Z}[\sqrt{D}]$ such that $px+zy=1$. Then,
	\begin{eqnarray*}
		m=pxm+zmy=pxm+ply=p(xm+ly)\in\langle p\rangle,
	\end{eqnarray*}
which  is a contradiction.  However, if $a$ is a unit, then $\langle  p\rangle=\langle ca\rangle=\langle c\rangle$. This means that $c=pn$ for some $n\in\mathbb{Z}[\sqrt{D}]$. Since $\langle  p,z\rangle=\langle c\rangle$, there exists $x,y\in\mathbb{Z}[\sqrt{D}]$ such that 
		$px+z=cy$. Then, we have 
		$$z=cy-px=pnb-px=p(nb-x)\in\langle p\rangle,$$
	which is a contradiction. Thus, $\langle p,z\rangle$ is a non-principal ideal.  
	
	On the other hand,  we suppose that  $\langle p,z\rangle$ is a non-principal ideal and assume for a contradiction that $z\notin I_p(D)$. 	By assumption, we have $zm\notin\langle p\rangle$, for any $m\in\mathbb{Z}[\sqrt{D}]$ such that $m\notin\langle p\rangle$. Since $z\notin\langle p\rangle$, $\bar{z}\notin\langle p\rangle$.  By choosing $m=\bar{z}$,  we conclude that $k:=zm=\lVert z\rVert\notin\langle p\rangle$; i.e., $k\in \mathbb{Z}$ and $\operatorname{gcd}(p,k)=1$. Then, by the Euclidean algorithm, there exists $x,y\in\mathbb{Z}$ such that $$1=px+ky=px+z\bar{z}y\in\langle p,z\rangle.$$ This implies that $\langle p,z\rangle$ is a principal ideal, which is a contradiction.  So, $z\in I_p(D)$.
\end{proof}

We further observe that:
\begin{prop}
	Let $z\in I_p(D)$.  Then,
\begin{enumerate}
	\item $p\mid\lVert z\rVert$ and 
	\item $\bar{z}\in I_p(D)$.
\end{enumerate}
Moreover, for any integer $t\in \mathbb{Z}$, $$  (\mathbb{Z}\cup R^*\cup tR^*)\cap I_p(D)=\emptyset,$$
where $R^*=\mathbb{Z}[\sqrt{D}]^{\times}$.
\end{prop}
\begin{proof}
		Let $z\in I_p(D)$. To prove (1), we assume for a contradiction that  $p\nmid\lVert z\rVert$.  Then $\operatorname{gcd}(p,\lVert z\rVert)=1$ and thus, by the Euclidean algorithm,  there exist $x,y\in\mathbb{Z}$ such that 
		$px+\lVert z\rVert y=1.$  This means that $		px+z\bar{z}y=1$. So,
		$\langle p,z\rangle=\langle 1\rangle$ is a principal ideal, which is a contradiction. 
		
		To prove (2), we again assume for a contradiction that  $\bar{z}\notin I_p(D)$. Then, by Proposition \ref{propnz}, $\langle p,\bar{z}\rangle$ is a principal ideal and hence so is $\langle p,z\rangle$ which leads to a contradiction.
		
		For the last assertion, we assume for a contradiction that there exists $a\in(\mathbb{Z}\cup R^*\cup tR^*)\cap I_p(D)$ for some  $t\in \mathbb{Z}$.   If $a\in\mathbb{Z}$, then $a^2=\lVert a\rVert$ and by (1), we have $p\mid a^2$.  So, $p\mid a$ and thus $a\in\langle p\rangle$.  This means that $a\notin I_p(D)$, a contradiction. If $a\in R^*$, then $\langle p, a\rangle$ is principal; i.e., $a\notin I_p(D)$,  a contradiction.  Finally, we suppose that $a=tu$ for some $u\in R^*$.  If $p\mid t$, then $\langle p, a\rangle=\langle p, t\rangle=\langle p\rangle$, a contradiction.  Also, if $p\nmid t$, then $p\nmid \lVert a \rVert$, since $$t^2=\lVert t \rVert=\lVert t \rVert\lVert u \rVert=\lVert tu \rVert=\lVert a \rVert.$$  So, by (1), $a\notin I_p(D)$, a contradiction.  
 \end{proof}

Note that $A(p,\bar{z})=A(p,z)^T$ and the idempotent factorization is preserved under the transpose operator, which leads us also conclude that:
\begin{cor}\label{cor23} For $z\in I_p(D)$, 
$A(p,z)$ satisfies conjecture \ref{con1} if and only if   $A(p,\bar{z})$ satisfies conjecture \ref{con1}.
\end{cor}
For $z\in I_p(D)$, denote $S_z:=\{m\in \mathbb{Z}[\sqrt{D}]\;|\; m\notin \langle p\rangle, zm\in \langle p\rangle \}$.  By Corollary \ref{cor23}, if $m=\bar{z}$ and $A(p,z)$ satisfies the conjecture, then  $A(p,\bar{z})$ must satisfies the conjecture.  However, if $m\in S_z$ with $m\neq \bar{z}$ , then the previous conclusion need not be true.  For example, we have $z:=1+\sqrt{10} \in I_3(10), m:=1+2\sqrt{10}\in S_z$ with $m\neq \bar{z}$.   By Theorem \ref{thmnz=p^2}, we can see that $A(3,z)$ satisfies the conjecture but, by Corollary \ref{coroppose}, $
A(3,m) $ does not satisfy the conjecture.

Moreover, by using the fact that idempotent factorization is also preserved under the similarity, it is enough to consider the factorization of  $A(p,z_0)$, where $z_0=z_1+z_2\sqrt{D}\in I_p(D)$ is a fundamental solution of the generalized Pell's equation $x^2-y^2D=\lVert z_0 \rVert$. 

\begin{prop} \label{propunit}
Let $z_0\in I_p(D)$ be a fundamental solution of a generalized Pell's equation $x^2-y^2D=\lVert z_0 \rVert$.  If $z\in I_p(D)$  is  in the same equivalence class of $z_0$, then,  \begin{center}
	$A(p,z_0)\in\operatorname{ID}_2(\mathbb{Z}[\sqrt{D}])$ if and only if  $A(p,z)\in\operatorname{ID}_2(\mathbb{Z}[\sqrt{D}])$.
\end{center} Furthermore, $A(p,z_0)$  satisfies conjecture \ref{con1} if and only if  $A(p,z)$  satisfies conjecture \ref{con1}.
\end{prop}
\begin{proof}
	By the assumption, $z=z_0u$ for some $u\in\mathbb{Z}[\sqrt{D}]$.  Then, by choosing the matrix $\begin{pmatrix}
	1&0\\0&u
	\end{pmatrix}\in M_2(\mathbb{Z}[\sqrt{D}])$ which is an invertible matrix, we see that
	\begin{eqnarray*}
		\begin{pmatrix}
			1&0\\0&\bar{u}
		\end{pmatrix}A(p,z_0)\begin{pmatrix}
			1&0\\0&u
		\end{pmatrix}=A(p,z_0u)=A(p,z);
	\end{eqnarray*}
namely, $A(p,z)$ and $A(p,z_0)$ are similar. Since  idempotent factorization is  preserved by  similarity, we reach to the conclusion. 
	Furthermore, if $A(p,z_0)$  satisfies conjecture \ref{con1}, then, there exists $a,b,c\in\mathbb{Z}[\sqrt{D}]$ with $a(1-a)=bc$ such that $$A(p,z_0)=\begin{pmatrix}
	a&b\\ c&1-a
	\end{pmatrix}\begin{pmatrix}
	\bar{a}&\bar{c}\\ \bar{b}& 1-\bar{a}
	\end{pmatrix}.$$ By the above similarity factorization, we have
	$$A(p,z)=	\begin{pmatrix}
	1&0\\0&\bar{u}
	\end{pmatrix}A(p,z_0)\begin{pmatrix}
	1&0\\0&u
	\end{pmatrix}=\begin{pmatrix}
	a&b\bar{u}\\cu&1-a
	\end{pmatrix}\begin{pmatrix}
	\bar{a}&\bar{c}\bar{u}\\ \bar{b}u& 1-\bar{a}
	\end{pmatrix},$$ which means that
	$A(p,z)$  satisfies the conjecture.   The converse can be proved by using the similar arguments.
\end{proof}

For example,  by Corollary \ref{corp=2}, we have $
A	(2,\sqrt{10})$ admits idempotent factorization.  Note that $19+6\sqrt{10}\in\mathbb{Z}[\sqrt{10}]^{\times}$. Then, by Proposition \ref{propunit}, we obtained that  $$
A(	2,\sqrt{10}\cdot(19+6\sqrt{10})^n)$$ admits idempotent factorization for any positive integer $n$.  Moreover, there is only one  fundamental solution class for  $x^2-10y^2=-10$ (see \cite{web1}).  Then, for any $z\in I_2(10)$ with $\lVert z \rVert=-10$, it turns out that $A(2,z)\in\operatorname{ID}_2(\mathbb{Z}[\sqrt{10}])$.

Recall that a matrix $M\in\mathbb{M}_2(R)$ is called \textbf{column--row} if there exists $a,b,c,d\in R$ such that 
$$M=\begin{pmatrix}
	a&0\\
	b&0
\end{pmatrix}\begin{pmatrix}
	c&d\\
	0&0
\end{pmatrix}.$$  Salce and Zanardo \cite{sz14}, in 2014, showed that if a singular matrix $M\in\mathbb{M}_2(R)$ with the ideal generated by its first row is principal, then $M$ is a column--row matrix (see Proposition 2.2 in \cite{sz14} for details).  This conclusion is extended to any row or any column in Corollary 4.1 in \cite{cz20}.    Unfortunately, the following proposition demonstrates that applying the similarity to a column-row matrix does not provide a new conclusion on idempotent factorization. 
\begin{prop}\label{crm} Let $A\in\mathbb{M}_2(\mathbb{Z}[\sqrt{D}])$ be a singular matrix.
	If $A$ is similar to a matrix having at least one row or column whose elements generate a principal ideal, then $A$ must be a column-row matrix.
\end{prop}
\begin{proof}
 Suppose that $A$ is similar to $\begin{pmatrix}
		s&0\\
		t&0
	\end{pmatrix}\begin{pmatrix}
		u&v\\
		0&0
	\end{pmatrix},$ for some $s,t,u,v\in \mathbb{Z}[\sqrt{D}]$.  So, there exists $\begin{pmatrix}
		a&b\\
		c&d
	\end{pmatrix}$ with $ad-bc=h\in \mathbb{Z}[\sqrt{D}]^{\times}$ for some $a,b,c,d\in \mathbb{Z}[\sqrt{D}]$ such that 
	\begin{eqnarray*}
		A&=&\begin{pmatrix}
			a&b\\
			c&d
		\end{pmatrix}\begin{pmatrix}
			s&0\\
			t&0
		\end{pmatrix}\begin{pmatrix}
			u&v\\
			0&0
		\end{pmatrix}\begin{pmatrix}
			dh^{-1}&-bh^{-1}\\
			-ch^{-1}&ah^{-1}
		\end{pmatrix}\\
		&=&\begin{pmatrix}
			as+bt&0\\
			cs+dt&0
		\end{pmatrix}\begin{pmatrix}
			udh^{-1}-vch^{-1}&-ubh^{-1}+vah^{-1}\\
			0&0
		\end{pmatrix},
	\end{eqnarray*}	
	which means that $A$ is a column-row matrix.
\end{proof}
 

\section{Main Results}

In this section, for any $z\in I_p(D)$, we are mainly interested in the conditions that $A(p,z)$ satisfies conjecture \ref{con1}. 
It is a direct calculation to see that $A(p,z) $ can be written as a product of two idempotent matrices, 
$$A(p,z)=\begin{pmatrix}
	a&b\\ c&1-a
\end{pmatrix}\begin{pmatrix}
	d&e\\ f& 1-d
\end{pmatrix},$$  if and only if,
\begin{eqnarray*}
	p&=&ad+bf\\
	z&=&ae+b(1-d)\\
	\bar{z}&=&cd+f(1-a)\\
	k:=\frac{\lVert z \rVert}{p}&=&ce+(1-a)(1-d)\\
	a(1-a)&=&bc\\
	d(1-d)&=&ef,
\end{eqnarray*} for some  $a,b,c,d,e,f\in\mathbb{Z}[\sqrt{D}]$.  Some of the above relations are redundant which can be reduced to:
\begin{lem} \label{lem01}
	Let $z\in I_p(D)$.  For any $a,b,c,d,e,f\in\mathbb{Z}[\sqrt{D}]$, we have that
	$$A(p,z)=BC \hbox{ where } B:=\begin{pmatrix}
		a&b\\ c&1-a
	\end{pmatrix} \hbox{ and } C:= \begin{pmatrix}
		d&e\\ f& 1-d
	\end{pmatrix},$$ are idempotent matrices  if and only if \begin{eqnarray}
		\label{000}	p&=&ad+bf\\
		\label{001}	z(1-a)&=&kb\\
		\label{002}	\bar{z}a&=&pc\\
		\label{003}	zd&=&pe\\
		\label{004}	\bar{z}(1-d)&=&kf.
	\end{eqnarray}
\end{lem}
\begin{proof}
	Since $B,C$ are idempotent,  $BA(p,z)=A(p,z)$ and $A(p,z)C=A(p,z)$.   For  $BA(p,z)=A(p,z)$  we have $(B-I)A(p,z)=0$.  Then,  \begin{eqnarray}
		(a-1)p&=&-b\bar{z}\label{01}\\
		(a-1)z&=&-bk\label{02}\\
		cp&=&a\bar{z}\label{03}\\
		cz&=&bk\label{04}\\
		a(1-a)&=&bc\label{05}.
	\end{eqnarray}
	We observe that  (\ref{01}) can be obtained from (\ref{02}) multiplied by $\bar{z}$ and (\ref{03}) can be obtained from (\ref{04}) multiplied by $\bar{z}$.  Also, (\ref{05}) comes from (\ref{02})  and  (\ref{03}).  Then, the above relations can be reduced to: 
	 \begin{eqnarray*}
		(1-a)z&=&kb\label{306}\\ \bar{z}a&=&pc\label{307}.
	\end{eqnarray*} Similarly, when we apply the same procedure to $A(p,z)(I-C) = 0$, we will get the relations:
	\begin{eqnarray*}
		(1-d)\bar{z}&=&kf\label{306}\\ zd&=&pe\label{307}.
	\end{eqnarray*}

	 For the converse, suppose that the relations (\ref{000}) to (\ref{004}) hold true.   Then, the relations (\ref{001}) and (\ref{002}) (multiplying them) yields that $B$ is an idempotent matrix.  Also, the relations (\ref{003}) and (\ref{004}) yields that $C$ is an idempotent matrix.   By multiplying equations (\ref{002}), (\ref{003}) and (\ref{004}) by $1-a,1-d$ and $e$, respectively, and use $p=ad+bf$ we obtain respectively that
	\begin{eqnarray*}
		\bar{z}&=&cd+(1-a)f\\
		z&=&ae+b(1-d)\\
		k&=&ce+(1-a)(1-d),
	\end{eqnarray*}
which completes the proof.
\end{proof}

Recall that the quadratic ring of integer $\mathbb{Z}[\sqrt{D}]$ of the quadratic number field  $\mathbb{Q}[\sqrt{D}]$ can be expressed explicitly depending on $D$.  Precisely,  (see, e.g., Theorem 2.2 in \cite{qnf}),
\begin{equation*}
	\mathbb{Z}[\sqrt{D}]	=	\begin{cases}
		\{a+b\sqrt{D}:a,b\in\mathbb{Z}\}&\text{if } D\equiv2,3(\operatorname{mod}4),\\
		\{\frac{a+b\sqrt{D}}{2}:a,b\in\mathbb{Z},a\equiv b(\operatorname{mod}2)\}&\text{if } D\equiv1(\operatorname{mod}4).
	\end{cases} 
\end{equation*}
So, we divide our calculation into 3 cases regarding to $D$;
	\begin{itemize}
	\item[case 1:] $D\equiv1(\operatorname{mod}4)$ with $z$ in the form $\frac{z_1+z_2\sqrt{D}}{2}$ and $z_1,z_2$ must be odd,
	\item[case 2:] $D\equiv1(\operatorname{mod}4)$ with $z$ in the form $z_1+z_2\sqrt{D}$,
	\item[case 3:] $D\equiv2,3(\operatorname{mod}4)$.
\end{itemize} 

Now, for a given $A(p,z)\in \mathbb{M}_2(\mathbb{Z}[\sqrt{D}])$  to find a necessary and sufficient conditions for the idempotent factorization into two idempotent matrices, by  Lemma \ref{lem01}, it suffices to verify the existence of $a,b,c\in\mathbb{Z}[\sqrt{D}]$ that correspond to the equations (\ref{001}) and (\ref{002}) and $d,e,f\in\mathbb{Z}[\sqrt{D}]$ that correspond to the equations (\ref{003}) and (\ref{004}) and also $a,d,b,f$ correspond to the equations (\ref{000}).

\textbf{Case 1:} The condition that $a_1,a_2,b_1,b_2,c_1,c_2\in\mathbb{Z}$ corresponding to \begin{eqnarray*}
	(\frac{z_1+z_2\sqrt{D}}{2})(1-\frac{a_1+a_2\sqrt{D}}{2})&=&k(\frac{b_1+b_2\sqrt{D}}{2})\\
	(\frac{z_1-z_2\sqrt{D}}{2})(\frac{a_1+a_2\sqrt{D}}{2})&=&p(\frac{c_1+c_2\sqrt{D}}{2}),
\end{eqnarray*}
is equivalent to
\begin{eqnarray*}
	\begin{matrix}
		z_1a_1-z_2Da_2&-2pc_1&\,&\,&\,&=&0\\
		-z_2a_1+z_1a_2&\,&-2pc_2&&\,&=&0\\
		z_1a_1+z_2a_2D& \,&\,&+2kb_1 &\,&=&2z_1\\
		z_2a_1+z_1a_2&\,&\,&\,&+2kb_2&=&2z_2.
	\end{matrix}
\end{eqnarray*}
To verify the existence of $a_1,a_2,c_1,c_2,b_1,b_2\in\mathbb{Z}$, we first solve the above system of linear equations (with variables $a_1,a_2,b_1,b_2,c_1,c_2$) by using the Gauss-Jordan elimination method (over $\mathbb{Q}[\sqrt{D}]$ ) and then concentrate only integral solutions.   By direct calculation, we have the row reduce echelon matrix for the above system as: $$\begin{pmatrix}
	1 & 0 & 0  &  0&\frac{z_1}{2p}&\frac{-z_2D}{2p} &\Big |&2\\
	0 & 1 & 0  &  0&\frac{-z_2}{2p}&\frac{z_1}{2p} &\Big |&0\\
	0 & 0 & 1  &  0&\frac{z_1^2+z_2^2D}{4p^2}&\frac{-2z_1z_2D}{4p^2} &\Big |&\frac{z_1}{p}\\
	0 & 0 & 0  &  1&\frac{-2z_1z_2}{4p^2}&\frac{z_1^2+z_2^2D}{4p^2} &\Big |&\frac{-z_2}{p}
\end{pmatrix}.$$ Then, for any $b_1,b_2\in\mathbb{Z}$, \begin{eqnarray}
	a_1&=&\frac{4p-z_1b_1+z_2b_2D}{2p}\label{a1/22}\\
	a_2&=&\frac{z_2b_1-z_1b_2}{2p}\label{a2/22}\\
	c_1&=&\frac{-(z_1^2+z_2^2D)b_1+2z_1z_2Db_2+4z_1p}{4p^2}\label{c1/22}\\
	c_2&=&\frac{2z_1z_2b_1-(z_1^2+z_2^2D)b_2-4z_2p}{4p^2}\label{c2/22}.\end{eqnarray}
Similarly, by repeating the above process with equations (\ref{003}) and (\ref{004}), for any $f_1,f_2\in\mathbb{Z}$, we obtained that 
\begin{eqnarray*}
	d_1&=&\frac{4p-z_1f_1-z_2f_2D}{2p}\label{d1/22}\\
	d_2&=&\frac{z_2f_1+z_1f_2}{2p}\label{d2/22}\\
	e_1&=&\frac{-(z_1^2+z_2^2D)b_1-2z_1z_2Df_2+4z_1p}{4p^2}\label{e1/22}\\
	e_2&=&\frac{2z_1z_2b_1+(z_1^2+z_2^2D)f_2-4z_2p}{4p^2}\label{e2/22}.\end{eqnarray*}
	By relation (\ref{000}), $a=(a_1+a_2\sqrt{D})/2, b=(b_1+b_2\sqrt{D})/2, d=(d_1+d_2\sqrt{D})/2$ and $f=(f_1+f_2\sqrt{D})/2$ must satisfy the condition $C_1$:
		\begin{eqnarray*}
			0&=&(p+k)b_1f_1+(p+k)Db_2f_2-z_1(b_1+f_1)+z_2(b_2-f_2)+4p-4p^2 \\
		0&=&4p^2(b_2f_1+b_1f_2)+z_1^2(b_2f_1+b_1f_2)-z_2^2D(b_2f_1+b_1f_2)-4pz_1(b_2+f_2)+4pz_2(b_1-f_1).
	\end{eqnarray*}
	
\textbf{Case 2:} 	The condition that $a_1,a_2,b_1,b_2,c_1,c_2\in\mathbb{Z}$ corresponding to \begin{eqnarray*}
	(z_1+z_2\sqrt{D})(1-\frac{a_1+a_2\sqrt{D}}{2})&=&k(\frac{b_1+b_2\sqrt{D}}{2})\\
	(z_1-z_2\sqrt{D})(\frac{a_1+a_2\sqrt{D}}{2})&=&p(\frac{c_1+c_2\sqrt{D}}{2}),
\end{eqnarray*}
is equivalent to
\begin{eqnarray*}
	\begin{matrix}
		
		z_1a_1-z_2Da_2&-pc_1&\,&\,&\,&=&0\\
		-z_2a_1+z_1a_2&\,&-pc_2&&\,&=&0\\
		z_1a_1+z_2a_2D& \,&\,&+kb_1 &\,&=&2z_1\\
		z_2a_1+z_1a_2&\,&\,&\,&+kb_2&=&2z_2.
	\end{matrix}
\end{eqnarray*}
   By direct calculation, we have the row reduce echelon matrix for the above system as: $$\begin{pmatrix}
	1 & 0 & 0  &  0&\frac{z_1}{p}&\frac{-z_2D}{p} &\Big |&2\\
	0 & 1 & 0  &  0&\frac{-z_2}{p}&\frac{z_1}{p} &\Big |&0\\
	0 & 0 & 1  &  0&\frac{z_1^2+z_2^2D}{p^2}&\frac{-2z_1z_2D}{p^2} &\Big |&\frac{2z_1}{p}\\
	0 & 0 & 0  &  1&\frac{-2z_1z_2}{p^2}&\frac{z_1^2+z_2^2D}{p^2} &\Big |&\frac{-2z_2}{p}
\end{pmatrix}.$$ Then, for any $b_1,b_2\in\mathbb{Z}$, \begin{eqnarray}
	a_1&=&\frac{2p-z_1b_1+z_2b_2D}{p}\label{a1/2}\\
	a_2&=&\frac{z_2b_1-z_1b_2}{p}\label{a2/2}\\
	c_1&=&\frac{-(z_1^2+z_2^2D)b_1+2z_1z_2Db_2+2z_1p}{p^2}\label{c1/2}\\
	c_2&=&\frac{2z_1z_2b_1-(z_1^2+z_2^2D)b_2-2z_2p}{p^2}\label{c2/2}.\end{eqnarray}
Similarly, by repeating the above process with equations (\ref{003}) and (\ref{004}), for any $f_1,f_2\in\mathbb{Z}$, we obtained that 
\begin{eqnarray*}
	d_1&=&\frac{2p-z_1f_1-z_2f_2D}{p}\label{d1/2}\\
	d_2&=&\frac{z_2f_1+z_1f_2}{p}\label{d2/2}\\
	e_1&=&\frac{-(z_1^2+z_2^2D)b_1-2z_1z_2Df_2+2z_1p}{p^2}\label{e1/2}\\
	e_2&=&\frac{2z_1z_2b_1+(z_1^2+z_2^2D)f_2-4z_2p}{p^2}\label{e2/2}.\end{eqnarray*}
	By relation (\ref{000}), $a=a_1+a_2\sqrt{D}, b=b_1+b_2\sqrt{D}, d=d_1+d_2\sqrt{D}$  and $f=f_1+f_2\sqrt{D}$ must satisfy the condition $C_2$:
		\begin{eqnarray*}
		0&=&(p+k)b_1f_1+(p+k)Db_2f_2-2z_1(b_1+f_1)+2z_2(b_2-f_2)+4p-4p^2\\
		0&=&p^2(b_2f_1+b_1f_2)+z_1^2(b_2f_1+b_1f_2)-z_2^2D(b_2f_1+b_1f_2)-2pz_1(b_2+f_2)+2pz_2(b_1-f_1).
	\end{eqnarray*}
	
	\textbf{Case 3:} 	The condition that $a_1,a_2,b_1,b_2,c_1,c_2\in\mathbb{Z}$ corresponding to \begin{eqnarray*}
	(z_1+z_2\sqrt{D})(1-a_1-a_2\sqrt{D})&=&k(b_1+b_2\sqrt{D})\\
	(z_1-z_2\sqrt{D})(a_1+a_2\sqrt{D})&=&p(c_1+c_2\sqrt{D}),
\end{eqnarray*}
is equivalent to
\begin{eqnarray*}
	\begin{matrix}
		
		z_1a_1-z_2Da_2&-pc_1&\,&\,&\,&=&0\\
		-z_2a_1+z_1a_2&\,&-pc_2&&\,&=&0\\
		z_1a_1+z_2a_2D& \,&\,&+kb_1 &\,&=&z_1\\
		z_2a_1+z_1a_2&\,&\,&\,&+kb_2&=&z_2.
	\end{matrix}
\end{eqnarray*} By direct calculation, we have the row reduce echelon matrix for the above system as:  $$\begin{pmatrix}
	1 & 0 & 0  &  0&\frac{z_1}{p}&\frac{-z_2D}{p} &\Big |&1\\
	0 & 1 & 0  &  0&\frac{-z_2}{p}&\frac{z_1}{p} &\Big |&0\\
	0 & 0 & 1  &  0&\frac{z_1^2+z_2^2D}{p^2}&\frac{-2z_1z_2D}{p^2} &\Big |&\frac{z_1}{p}\\
	0 & 0 & 0  &  1&\frac{-2z_1z_2}{p^2}&\frac{z_1^2+z_2^2D}{p^2} &\Big |&\frac{-z_2}{p}
\end{pmatrix}.$$ Then, for any $b_1,b_2\in\mathbb{Z}$, \begin{eqnarray}
	a_1&=&\frac{p-z_1b_1+z_2b_2D}{p}\label{a1}\\
	a_2&=&\frac{z_2b_1-z_1b_2}{p}\label{a2}\\
	c_1&=&\frac{-(z_1^2+z_2^2D)b_1+2z_1z_2Db_2+z_1p}{p^2}\label{c1}\\
	c_2&=&\frac{2z_1z_2b_1-(z_1^2+z_2^2D)b_2-z_2p}{p^2}\label{c2}.
\end{eqnarray} Similarly, by repeating the above process with equations (\ref{003}) and (\ref{004}), for any $f_1,f_2\in\mathbb{Z}$, we obtained that 
\begin{eqnarray*}
	\label{d1}	d_1&=&\frac{p-z_1f_1-z_2f_2D}{p}\\
	\label{d2}	d_2&=&\frac{-z_2f_1-z_1f_2}{p}\\
	\label{e1}	e_1&=&\frac{-(z_1^2+z_2^2D)f_1-2z_1z_2Df_2+z_1p}{p^2}\\
	\label{e2}	e_2&=&\frac{-2z_1z_2f_1-(z_1^2+z_2^2D)f_2+z_2p}{p^2}.
\end{eqnarray*}
	By relation (\ref{000}), $a=a_1+a_2\sqrt{D}, b=b_1+b_2\sqrt{D}, d=d_1+d_2\sqrt{D}$  and $f=f_1+f_2\sqrt{D}$ must satisfy the condition $C_3$:
		\begin{eqnarray*}
		0&=&(p+k)b_1f_1+(p+k)Db_2f_2-z_1(b_1+f_1)+z_2(b_2-f_2)+p-p^2\\
		0&=&p^2(b_2f_1+b_1f_2)+z_1^2(b_2f_1+b_1f_2)-z_2^2D(b_2f_1+b_1f_2)-pz_1(b_2+f_2)+pz_2(b_1-f_1).
	\end{eqnarray*}

From the above 3 cases, in particular, if $f_1=b_1\text{ and }f_2=-b_2$, then $a=\bar{d}\text{ and }b=\bar{f}$.   This implies that the conditions $C_1$, $C_2$ and $C_3$ become: 

	\begin{eqnarray}\label{last1}
		0&=&(p+k)b_1^2-(p+k)Db_2^2-2z_1b_1+2z_2b_2D+4p-4p^2,
	\end{eqnarray}

	\begin{eqnarray}\label{last2}
		0&=&(p+k)b_1^2-(p+k)Db_2^2-4z_1b_1+4z_2b_2D+4p-4p^2,
	\end{eqnarray}
and
	\begin{eqnarray}\label{last3}
		0&=&(p+k)b_1^2-(p+k)Db_2^2-2z_1b_1+2z_2Db_2+p-p^2,
	\end{eqnarray}
respectively, where $k=\lVert z \rVert/p$.  We now reach to a necessary and sufficient condition for the idempotent factorization in conjecture \ref{con1} in term of a system of quadratic Diophantine equations (instead of a system of equations in the ring of integers).

\begin{thm}\label{mainthm} Let $z\in I_p(D)$. Then, $A(p,z)$ satisfies  conjecture \ref{con1} if and only if  
\begin{itemize}
	\item[case 1:] the system of equations (\ref{a1/22}) to  (\ref{c2/22}) and (\ref{last1}) has an integral solution, 
		\item[case 2:] the system of equations (\ref{a1/2}) to  (\ref{c2/2}) and (\ref{last2}) has an integral solution, 
\item[case 3:] the system of equations (\ref{a1}) to  (\ref{c2}) and (\ref{last3}) has an integral solution.		
	\end{itemize}
\end{thm}

Note from the above theorem that, in general, the condition for the idempotent factorization depends on the fives equations.  However, in particular $p, D$ and $z\in I_p(D)$, only one equation  is sufficient to determine the idempotent factorization. 

\begin{cor}\label{corp=2} Let $z\in I_2(D)$ with
	$D\equiv2(\operatorname{mod4})$ and $\lVert z\rVert=2k$ for some $k\in\mathbb{Z}$ where $k\equiv3(\operatorname{mod4})$, then  $A(2,z)$ satisfies  conjecture \ref{con1} if and only if the quadratic Diophantine equation \begin{equation}
		0=(2+k)b_1^2-(2+k)Db_2^2-2z_1b_1+2z_2Db_2-2\label{p=2}
	\end{equation} has an integer solution.
\end{cor}
\begin{proof} Suppose that $A(2,z)$ satisfies the conjecture. By Theorem \ref{mainthm} case 3 with $p=2$, we conclude  from (\ref{last3}) that the equation (\ref{p=2}) has an integral solution. 
	Conversely, we assume that the equation (\ref{p=2}) has a solution. Since $k\equiv3(\operatorname{mod4})$, we have 
	$\lVert z\rVert\equiv2(\operatorname{mod4})$.   Let $z=z_1+z_2\sqrt{D}$.  So, $z_1^2+2z_2^2\equiv2(\operatorname{mod4})$, which implies that $z_1\equiv0(\operatorname{mod2})$ and $z_2\equiv1(\operatorname{mod2})$.    Since the equation (\ref{p=2}) has a solution $b_1,b_2\in\mathbb{Z}$, $0\equiv b_1^2+2b_2^2+2(\operatorname{mod4})$. This yields that $b_1\equiv0(\operatorname{mod2})$ and $b_2\equiv1(\operatorname{mod2})$.
	Now, we can write $z_1=2h_1,z_2=2h_2+1,b_1=2h_3,b_2=2h_4+1$ and $D=4h_5+2$ for some $h_1,h_2,h_3,h_4,h_5\in\mathbb{Z}$. By equations (\ref{a1}) to (\ref{c2}), we have
	\begin{eqnarray*}
		a_1&=&1-2h_1h_3+(2h_5+1)(2h_4+1)(2h_2+1)\\
		a_2&=&h_3(2h_2+1)-h_1(2h_4+1)\\
		c_1&=&-h_3(2h_1^2+(2h_2+1)^2(2h_5+1))+2h_1(2h_2+1)(2h_4+1)(2h_5+1)+h_1 \\
		c_2&=&2h_1h_3(2h_2+1)-h_1^2(2h_4+1)-(2h_2+1)(2h_2h_5+h_2+h_5+1),
	\end{eqnarray*} are all integers.
\end{proof}

To show that conjecture \ref{con1} holds true when $\lVert z\rVert=-p^2$, the following lemmas are required.
\begin{lem}\label{lemgcd1}
	Let $z\in I_p(D)$ with $\lVert z\rVert=-p^2$.  Then,  $\operatorname{gcd}(z_1,z_2D)=1,$ for $z=z_1+z_2\sqrt{D}$ or  $z=(z_1+z_2\sqrt{D})/2$.
\end{lem}
\begin{proof}
We first consider $z\in I_p(D)$ in the form $z_1+z_2\sqrt{D}.$  Let $d=\operatorname{gcd}(z_1,z_2D)$. Then, $d\mid z_1$ and $d\mid z_2D.$ This implies that $$d\mid z_1^2-Dz_2^2=-p^2.$$ So, the possibility of $d$ are $1,p,p^2$. If $d=p$, then $p\mid z_1$ and $p\mid z_2D.$ This implies that $p^2\mid Dz_2^2$, since $p^2\mid z_1^2$ and $p^2\mid z_1^2-Dz_2^2$. If $p\nmid D$, then $p^2\mid z_2^2$ and thus $p\mid z_2$; i.e., $p\mid z$, a contradiction. However, if $p\mid D$, then $p\mid z_2^2$, which again leads to a contradiction.   If  $d=p^2$, we have $p^2\mid z_1$ and $p^2\mid z_2D$. Since $D$ is square free, $p\mid z_1$ and $p\mid z_2,$ a contradiction. It follows that $d=1.$  

In the other case, when $D\equiv1(\operatorname{mod}4)$ and $z$ is in the form $\frac{z_1+z_2\sqrt{D}}{2}$ with $z_1,z_2\equiv1(\operatorname{mod}2)$, we have that $z_1^2-z_2^2D=-4p^2$.  Now, we obtained that $d\mid z_1^2-Dz_2^2=-4p^2,$ where  $d=\operatorname{gcd}(z_1,z_2D)$. So, the possibility of $d$ are $1,2,p,2p,4p,p^2, 2p^2, 4p^2$.  If $2\mid d$, then $2\mid z_1$, a contradiction to the fact that $z_1\equiv1(\operatorname{mod}2)$. If $p\mid d$, then $p\mid z_1$ and $p\mid z_2D.$ This implies that $p\mid z_2$, since $p^2\mid z_1^2$ and $p^2\mid z_1^2-Dz_2^2$ and $D$ is square free, a contradiction. This implies that $d=1.$
\end{proof}

The next lemma shows that if $z\in I_2(D)$, then $\lVert z \rVert\neq-4$.
\begin{lem}\label{lem02}
		Let $z\in \mathbb{Z}[\sqrt{D}]$ with $\lVert z \rVert=-4$. Then, $\langle2,z\rangle$ is a principal ideal.
\end{lem}
\begin{proof}By assumption we have $z_1^2-Dz_2^2=-4$.
We first consider the case where $D$ is odd square free integer. \begin{itemize}
		\item [Case1:] $z_1$ is even. Then,  $z_2$  must even, which implies that $z\in\langle 2\rangle$.
		\item [Case2:] $z_1$ is odd. Then,  $z_2$ must be odd. Then, there exists $m,n\in\mathbb{Z}$ such that 
		\begin{eqnarray*}
			-4&=&(2m+1)^2-D(2n+1)^2\\
			&=&4m^2+4m-1+D(-4n^2-4n-1)\\
			&=&4(m^2-Dn^2+m-nD)+1-D.
		\end{eqnarray*}
So, $D\equiv1(\operatorname{mod}4)$.  Now, we obtained that
		\begin{eqnarray*}
		z=	z_1+z_2\sqrt{D}&=&(2m+1)+(2n+1)\sqrt{D}\\
						&=&2(m+n\sqrt{D}+(\frac{1+\sqrt{D}}{2}))\\
			&\in&\langle 2\rangle.
	\end{eqnarray*}\end{itemize} 
	For the case where $D$ is even square free integer.  By the assumption that $z_1^2-Dz_2^2=-4$,  $z_1$ must be even.  If $z_2$ is also even, we complete the proof. If $z_2$ is odd, then there exists $m,n\in\mathbb{Z}$ such that 
	\begin{eqnarray*}
		-4&=&4m^2-D(4n^2+4n+1)\\
		&=&4(m^2-Dn^2-Dn)-D.
	\end{eqnarray*}
	It follows that $D\equiv0(\operatorname{mod}4)$, which is a contradiction. 
\end{proof}
By Lemma \ref{lem02}, the prime numbers in the following theorem are odd.
\begin{thm}\label{thmnz=p^2}
Let $z\in I_p(D)$  with  $\lVert z\rVert=-p^2$. Then, \begin{center}
	 $A(p,z)$ satisfies conjecture \ref{con1}.
\end{center}
\end{thm}
\begin{proof}
	In this proof, we divided the proof into two parts according to the form of $z$: $z=(z_1+z_2\sqrt{D})/2$ and $z=z_1+z_2\sqrt{D}$ where $z_1,z_2\in \mathbb{Z}$.  For the form $z=\frac{z_1+z_2\sqrt{D}}{2}$, by Theorem \ref{mainthm} with $k=\lVert z\rVert/p=-p$, the equation (\ref{last1}) becomes \begin{eqnarray*}
		2p-2p^2&=&z_1b_1-z_2b_2D.
	\end{eqnarray*}
	By Lemma \ref{lemgcd1} and the Euclidean algorithm, there exists $x,y\in\mathbb{Z}$ such that $1=z_1x+(-z_2D)y.$ Now, for any $m\in\mathbb{Z}$, we let,  \begin{eqnarray*}
		b_1&=&2x(p-p^2)-z_2mD\\
		b_2&=&2y(p-p^2)-mz_1.
	\end{eqnarray*}
Since $D\equiv 1(\operatorname{mod}4)$ and $z_1\equiv z_2(\operatorname{mod}2)$, $b_1\equiv b_2(\operatorname{mod}2)$.   By substituting $b_1,b_2$ in equation (\ref{a1/22}), we obtained that
	\begin{eqnarray*}
		a_1&=&\frac{4p-z_1(2x(p-p^2)-z_2mD)+z_2(2y(p-p^2)-mz_1)D}{2p}\\
		&=&\frac{4p-z_1(2x(p-p^2)+z_2(2y(p-p^2))D}{2p}\\
		&=&2(1-\frac{1-p}{2}z_1x+\frac{1-p}{2}z_2Dy)\in\mathbb{Z}.\\
	\end{eqnarray*}
	Similarly, by substituting $b_1,b_2$ in equation (\ref{a2/22}), we have that \begin{eqnarray*}
		a_2&=&\frac{z_2(2x(p-p^2)-z_2mD)-z_1(2y(p-p^2)-mz_1)}{2p}\\
		&=&\frac{-4mp^2+z_2(2x(p-p^2))-z_1(2y(p-p^2))}{2p}\\
		&=&-2mp+z_2x(1-p))-z_1(y(1-p))\\
		&=&2(-mp+z_2x\frac{1-p}{2}-z_1y\frac{1-p}{2})\in\mathbb{Z}.
	\end{eqnarray*} 
	Since $p$ is odd, $(1-p)/2\in \mathbb{Z}$ and thus  $a_1\equiv a_2(\operatorname{mod}2)$.  Note that $z_1^2+z_2^2D=2z_1^2+4p^2$, because $\lVert z \rVert=(z_1^2-z_2^2D)/4=-p^2$.  Then,  by substituting $b_1,b_2$ in equation (\ref{c1/22}), we have that
	\begin{eqnarray*}
		c_1&=&\frac{-(z_1^2+z_2^2D)(2x(p-p^2)-z_2mD)+2z_1z_2D(2y(p-p^2)-mz_1))+4pz_1}{4p^2}\\
		&=&\frac{-(2z_1^2+4p^2)(2x(p-p^2)-z_2mD)+2z_1z_2D(2y(p-p^2)-mz_1))+4pz_1}{4p^2}\\
		&=&\frac{-4z_1^2x(p-p^2)+2z_1^2z_2mD-8p^2x(p-p^2)+4p^2z_2mD+4z_1z_2Dy(p-p^2)-2z_1^2z_2mD+4pz_1}{4p^2}\\
		&=&\frac{-4z_1^2x(p-p^2)-8p^2x(p-p^2)+4p^2z_2mD+4z_1z_2Dy(p-p^2)+4pz_1}{4p^2}\\
		&=&\frac{4pz_1(-z_1x(1-p)+z_2Dy(1-p)+1)-8p^2x(p-p^2)+4p^2z_2mD}{4p^2}\\
		&=&\frac{4p^2z_1-8p^2x(p-p^2)+4p^2z_2mD}{4p^2}\\
		&=&z_1-2x(p-p^2)+z_2mD\in\mathbb{Z}.
	\end{eqnarray*}
	Note that $z_1^2+z_2^D=-4p^2+2z_2^2D$, because $\lVert z \rVert=-p^2$. Then, by substituting $b_1,b_2$ in equation (\ref{c2/22}), we have that
	\begin{eqnarray*}
		c_2&=&\frac{2z_1z_2(2x(p-p^2)-z_2mD)-(z_1^2+z_2^2D)(2y(p-p^2)-mz_1))-4pz_2}{4p^2}\\
		&=&\frac{2z_1z_2(2x(p-p^2)-z_2mD)-(-4p^2+2z_2^2D)(2y(p-p^2)-mz_1))-4pz_2}{4p^2}\\
		&=&\frac{4z_1z_2x(p-p^2)-2z_1z_2^2mD+8p^2y(p-p^2)-4p^2mz_1-4z_2^2yD(p-p^2)+2z_1z_2^2Dm-4pz_2}{4p^2}\\
		&=&\frac{4z_1z_2x(p-p^2)+8p^2y(p-p^2)-4p^2mz_1-4z_2^2yD(p-p^2)-4pz_2}{4p^2}\\
		&=&\frac{4pz_2(z_1x(1-p)-z_2Dy(1-p)-1)+8p^2y(p-p^2)-4p^2mz_1}{4p^2}\\
		&=&\frac{-4p^2z_2+8p^2y(p-p^2)-4p^2mz_1}{4p^2}\\
		&=&-z_2+2y(p-p^2)-mz_1\in\mathbb{Z}.
	\end{eqnarray*}
	Since $z_1\equiv z_2(\operatorname{mod}2)$, we have that $c_1\equiv c_2(\operatorname{mod}2)$, for any $m\in\mathbb{Z}$. By Theorem \ref{mainthm}, we have $A(p,z)$ admits the idempotent factorization.

For the form $z=z_1+z_2\sqrt{D}$ with $D\equiv1(\operatorname{mod4})$ and $\lVert z\rVert=-p^2$, the equation (\ref{last2}) becomes
\begin{eqnarray*}
	p-p^2=z_1b_1-z_2Db_2.
\end{eqnarray*}
By Lemma \ref{lemgcd1} and the Euclidean algorithm, there exists $x,y\in\mathbb{Z}$ such that $1=z_1x+(-z_2D)y$.  Now, for any $m\in\mathbb{Z}$, we let,  \begin{eqnarray*}
	b_1&=&x(p-p^2)-z_2mD\\
	b_2&=&y(p-p^2)-mz_1.
\end{eqnarray*}
Since $D\equiv 1(\operatorname{mod}4)$ and $z_1\equiv z_2(\operatorname{mod}2)$, $b_1\equiv b_2(\operatorname{mod}2)$.   By substituting $b_1,b_2$ in equation (\ref{a1/2}), we obtained that
\begin{eqnarray*}
	a_1&=&\frac{2p-z_1(x(p-p^2)-z_2mD)+z_2(y(p-p^2)-mz_1)D}{p}\\
	&=&\frac{2p-z_1x(p-p^2)+z_2(y(p-p^2))D}{p}\\
	&=&2(1-\frac{1-p}{2}z_1x+\frac{1-p}{2}z_2Dy)\in\mathbb{Z}.\\
\end{eqnarray*}
Similarly, by substituting $b_1,b_2$ in equation (\ref{a2/2}), we have that \begin{eqnarray*}
	a_2&=&\frac{z_2(x(p-p^2)-z_2mD)-z_1(y(p-p^2)-mz_1)}{p}\\
	&=&\frac{-mp^2+z_2(x(p-p^2))-z_1(y(p-p^2))}{p}\\
	&=&-mp+z_2x(1-p))-z_1(y(1-p))\\
	&=&2(-mp+z_2x\frac{1-p}{2}-z_1y\frac{1-p}{2})\in\mathbb{Z}.
\end{eqnarray*} 
Since $p$ is odd, $(1-p)/2\in \mathbb{Z}$ and thus  $a_1\equiv a_2(\operatorname{mod}2)$.  Note that $z_1^2+z_2^2D=2z_1^2+p^2$, because $\lVert z \rVert=(z_1^2-z_2^2D)=-p^2$.  Then,  by substituting $b_1,b_2$ in equation (\ref{c1/2}), we have that
\begin{eqnarray*}
	c_1&=&\frac{-(z_1^2+z_2^2D)(x(p-p^2)-z_2mD)+2z_1z_2D(y(p-p^2)-mz_1))+2pz_1}{p^2}\\
	&=&\frac{-(2z_1^2+p^2)(x(p-p^2)-z_2mD)+2z_1z_2D(y(p-p^2)-mz_1))+2pz_1}{p^2}\\
	&=&\frac{-2z_1^2x(p-p^2)-p^2x(p-p^2)+p^2z_2mD+2z_1z_2Dy(p-p^2)+2pz_1}{p^2}\\
	&=&\frac{-2pz_1(p-p^2)(-z_1x+z_2Dy)-p^2x(p-p^2)+p^2z_2mD+2pz_1}{p^2}\\
	&=&\frac{2p^2z_1-p^2x(p-p^2)+p^2z_2mD}{p^2}\\
	&=&2z_1-2xp\left( \frac{1-p}{2} \right)+z_2mD\in\mathbb{Z}.
\end{eqnarray*}
Note that $z_1^2+z_2^D=-p^2+2z_2^2D$, because $\lVert z \rVert=-p^2$. Then, by substituting $b_1,b_2$ in equation (\ref{c2/2}), we have that
\begin{eqnarray*}
	c_2&=&\frac{2z_1z_2(x(p-p^2)-z_2mD)-(z_1^2+z_2^2D)(y(p-p^2)-mz_1))-2pz_2}{p^2}\\
	&=&\frac{2z_1z_2(x(p-p^2)-z_2mD)-(-p^2+2z_2^2D)(y(p-p^2)-mz_1))-2pz_2}{p^2}\\
	&=&\frac{2z_1z_2x(p-p^2)+p^2y(p-p^2)-p^2z_1m-2z_2^2Dy(p-p^2)-2pz_1}{p^2}\\
	&=&\frac{2z_2(p-p^2)(z_1x-z_2Dy)+p^2y(p-p^2)-p^2z_1m-2pz_1}{4p^2}\\
	&=&\frac{-2p^2z_2+p^2y(p-p^2)-p^2z_1m}{4p^2}\\
	&=&-2z_2+2yp\left( \frac{1-p}{2} \right)-mz_1\in\mathbb{Z}.
\end{eqnarray*}
Since $z_1\equiv z_2(\operatorname{mod}2)$, we have  $c_1\equiv c_2(\operatorname{mod}2)$, for any $m\in\mathbb{Z}$. By Theorem \ref{mainthm}, we have that $A(p,z)$ admits the idempotent factorization.

	For the form $z=z_1+z_2\sqrt{D}$ with $D\equiv2,3(\operatorname{mod4})$ and $\lVert z\rVert=-p^2$, by Theorem \ref{mainthm} with $k=-p$, the equation (\ref{last3}) becomes \begin{eqnarray*}
		\frac{-p(p-1)}{2}&=&z_1b_1-z_2Db_2.
	\end{eqnarray*}
	Again, by Lemma \ref{lemgcd1} and the Euclidean algorithm, there exists $x,y\in\mathbb{Z}$ such that $1=z_1x+(-z_2D)y.$  Now, for any $m\in\mathbb{Z}$, we let, $m\in\mathbb{Z}$, \begin{eqnarray*}
		b_1&=&\frac{xp(1-p)}{2}-z_2mD\\
		b_2&=&\frac{yp(1-p)}{2}-mz_1.
	\end{eqnarray*}
	By substituting $b_1,b_2$ in equation (\ref{a1}), we obtained that
	\begin{eqnarray*}
		pa_1&=&p-z_1(\frac{xp(1-p)}{2}-z_2mD)+z_2(\frac{yp(1-p)}{2}-mz_1)D\\
	pa_1	&=&p-z_1\frac{xp(1-p)}{2}+z_2\frac{yp(1-p)}{2}\\
	pa_1	&=&p+\frac{p(1-p)}{2}\\
a_1	&=&1+\frac{p-1}{2}\in\mathbb{Z}.
	\end{eqnarray*}
	Similarly, by substituting $b_1,b_2$ in equation (\ref{a2}), we obtained that \begin{eqnarray*}
		pa_2&=&z_2(\frac{xp(1-p)}{2}-z_2mD)-z_1(\frac{yp(1-p)}{2}-mz_1)\\
	pa_2	&=&pz_2x\frac{1-p}{2}-pz_1y\frac{1-p}{2}+m(z_1^2-z_2^2D)\\
	pa_2	&=&pz_2x\frac{1-p}{2}-pz_1y\frac{1-p}{2}-mp^2\\
	a_2	&=&z_2x\frac{1-p}{2}-z_1y\frac{1-p}{2}-mp\in\mathbb{Z}.
	\end{eqnarray*} Denote that $z_1^2+z_2^2D=2z_1^2+p^2$, because $\lVert z \rVert = z_1^2-z_2^2D=-p^2$. Then,  by substituting $b_1,b_2$ in equation (\ref{c1}), we obtained that
	\begin{eqnarray*}
		p^2c_1&=&-(z_1^2+z_2^2D)(\frac{xp(1-p)}{2}-z_2mD)+2z_1z_2D(\frac{yp(1-p)}{2}-mz_1)+z_1p\\
		p^2c_1	&=&-(2z_1^2+p^2)(\frac{xp(1-p)}{2}-z_2mD)+2z_1z_2D(\frac{yp(1-p)}{2}-mz_1)+z_1p\\
		p^2c_1	&=&-2z_1^2x\frac{p(1-p)}{2}-p^2x\frac{p(1-p)}{2}+mp^2z_2D+2z_1z_2Dy\frac{p(1-p)}{2}+z_1p\\
	p^2c_1		&=&-2z_1(z_1x\frac{p(1-p)}{2}-z_2Dy\frac{p(1-p)}{2})+z_1p+p^2(-x\frac{p(1-p)}{2}+mz_2D)\\
	p^2c_1		&=&z_1p(p-1)+z_1p+p^2(-x\frac{p(1-p)}{2}+mz_2D)\\
	p^2c_1		&=&p^2(z_1-x\frac{p(1-p)}{2}+mz_2D)\\
	c_1	&=&z_1-x\frac{p(1-p)}{2}+mz_2D\in\mathbb{Z}.
	\end{eqnarray*}
	Denote that $z_1^2+z_2^2D=-p^2+2z_2^2D$, because $\lVert z \rVert=-p^2$. Then,  by substituting $b_1,b_2$ in equation (\ref{c2}), we have that
	\begin{eqnarray*}
		p^2c_2&=&2z_1z_2(\frac{xp(1-p)}{2}-z_2mD)-(z_1^2+z_2^2D)(\frac{yp(1-p)}{2}-mz_1)-z_2p\\
		p^2c_2	&=&2z_1z_2x\frac{p(1-p)}{2}+p^2(y\frac{p(1-p)}{2}-mz_1)-2z_2^2Dy\frac{p(1-p)}{2}-pz_2\\
		p^2c_2	&=&p^2(y\frac{p(1-p)}{2}-mz_1)+2z_2(z_1x\frac{p(1-p)}{2}-z_2Dy\frac{p(1-p)}{2})-pz_2\\
		p^2c_2	&=&p^2(y\frac{p(1-p)}{2}-mz_1)+2z_2(z_1x\frac{p(1-p)}{2}+z_2p(1-p)-pz_2\\
		p^2c_2	&=&p^2(y\frac{p(1-p)}{2}-mz_1)+2z_2(z_1x\frac{p(1-p)}{2}-z_2p^2\\
	c_2	&=&y\frac{p(1-p)}{2}-mz_1-z_2\in\mathbb{Z}.
	\end{eqnarray*}
By Theorem \ref{mainthm}, we can conclude that $A(p,z)$ admits idempotent factorization. 
\end{proof}

For example, $1+\sqrt{10}\in I_3(10)$,  Cossu and Zanardo \cite{cz20} expressed the idempotent factorization as: $$	A	(3,1+\sqrt{10})=\begin{pmatrix}
		2+2\sqrt{10}&7+\sqrt{10}\\-6&-1-2\sqrt{10}
	\end{pmatrix}\begin{pmatrix}
		2-2\sqrt{10}&-6
		\\7-\sqrt{10}&-1+2\sqrt{10}
	\end{pmatrix}.$$ By using the proof of Theorem \ref{thmnz=p^2} (choosing $m$ differently),  $A	(3,1+\sqrt{10})$ can also be factored as,  for example, 
	\begin{eqnarray*}
	A(3,1+\sqrt{10})
		&=&\begin{pmatrix}
			2+5\sqrt{10}&17+2\sqrt{10}\\ -16+\sqrt{10}&-1-5\sqrt{10}
		\end{pmatrix}\begin{pmatrix}
			2-5\sqrt{10}&-16-\sqrt{10}\\ 17-2\sqrt{10}& -1+5\sqrt{10}
		\end{pmatrix} \\
		&=&\begin{pmatrix}
			2-19\sqrt{10}&-36-6\sqrt{10}\\ 64-7\sqrt{10}&-1+19\sqrt{10}
		\end{pmatrix}\begin{pmatrix}
			2+19\sqrt{10}&64+7\sqrt{10}\\ -36+6\sqrt{10}& -1-19\sqrt{10}
		\end{pmatrix}\\
		&=&\begin{pmatrix}
			2+8\sqrt{10}&27+3\sqrt{10}\\ -26+2\sqrt{10}&-1-8\sqrt{10}
		\end{pmatrix}\begin{pmatrix}
			2-8\sqrt{10}&-26-2\sqrt{10}\\ 27-3\sqrt{10}& -1+8\sqrt{10}
		\end{pmatrix}\\
		&=&\begin{pmatrix}
			2-4\sqrt{10}&-13-1\sqrt{10}\\ 14-2\sqrt{10}&-1+4\sqrt{10}
		\end{pmatrix}\begin{pmatrix}
			2-4\sqrt{10}&14+2\sqrt{10}\\ -13+1\sqrt{10}& -1+4\sqrt{10}
		\end{pmatrix}, 
	\end{eqnarray*}
 etcetera.   From the above example, we can see that the idempotent factorization in the Cossu and Zanardo's conjecture need not be unique.

 In 2021, Matthews and Robertson \cite{mr21} gave a new method to solve a binary quadratic Diophantine equation $ax^2+bxy+cy^2+dx+ey+f=0$ by transforming it into the form $$as_2^2X^2+br_2s_2XY+cr_2^2Y^2=M,$$ where $M=-r_2^2s_2^2\left(ae^2-bde+cd^2+f\Delta_1\right)/\Delta_1,$ $r_1/r_2=\alpha/\Delta_1$ and $s_1/s_2=\beta/\Delta_1$ , so that $\operatorname{gcd}(r_1,r_2)=\operatorname{gcd}(s_1,s_2)=1$, when
$\Delta_1=b^2-4ac,$ $\alpha=2cd-be$ and $\beta=2ae-db$, the Florida transform.  So, by applying the Florida transform to equation (\ref{last3}), we obtain that, for $z\in I_p(D)$, 
\begin{eqnarray}
	X^2-DY^2=(p+k-1)p^2,\label{florida}
\end{eqnarray}
where $\operatorname{gcd}(p+k,z_1)=\operatorname{gcd}(p+k,z_2)=1$. We recall that a necessary condition for the solubility of $x^2-Dy^2=N$ is that $u^2\equiv D(\operatorname{mod}\vert N\vert)$ shall be soluble.  We also recall that, in number theory, the Kronecker symbol is a generalization of the Legendre symbol which can be expressed as:
\begin{eqnarray*}
	\left( \frac{a}{b} \right)=\begin{cases}
		1&\text{if $a$ is congruent to a perfect square modulo $b$}\\
		-1&\text{if $a$ is not congruent to any perfect square modulo $b$}\\
		0&\text{if $a\equiv0(\operatorname{mod }b)$.}
	\end{cases}
\end{eqnarray*}  From this observation, we find a necessary condition for the existence of  the  idempotent factorization that satisfy the conjecture. 

\begin{cor}\label{corflorida}
Let $z\in I_p(D)$ with $D\equiv2,3(\operatorname{mod4})$. Suppose that $p\nmid D$ and $\operatorname{gcd}(p+k,z_1)=\operatorname{gcd}(p+k,z_2)=1$. If the Kronecker symbol of $\left( \frac{D}{\lvert p+k-1\rvert} \right)=-1$, then  $A(p,z)$ does not satisfy conjecture \ref{con1}. 
\end{cor}
\begin{proof}
	If $\left( \frac{D}{\lvert p+k-1\rvert} \right)=-1$, then $x^2\equiv D (\operatorname{mod }\lvert p+k-1\rvert)$ has no solution. Then, the equation (\ref{florida}) has also no solution. This implies that the equation (\ref{last3}) has no solution. By Theorem \ref{mainthm}, $A(p,z)$ does not satisfy the  conjecture.
\end{proof}

Note that the Kronecker symbol of $\left( \frac{10}{\lvert 7\rvert} \right)=-1$. By Corollary \ref{corflorida}, we can conclude that  $A(3,5+\sqrt{10})$ and $A(3,-5+\sqrt{10})$  do not satisfy the  conjecture.  

Moreover, for $z\in I_p(D)$ with $D\equiv2(\operatorname{mod}4)$, we give an example that  $A(p,z)$ also oppose the  conjecture.
\begin{cor}\label{coroppose}
	Let $z\in I_p(D)$ with $D\equiv2(\operatorname{mod}4)$, $p\equiv3(\operatorname{mod}4)$and  $p+\sfrac{\lVert z\rVert}{p}\equiv2(\operatorname{mod}4).$ Then,  $A(p,z)$ does not satisfy conjecture \ref{con1}. 
\end{cor}
\begin{proof}
	By assumption, we have $p-p^2\equiv2(\operatorname{mod}4)$ and $\lVert z\rVert\equiv1(\operatorname{mod}4)$. Then, $p-p^2=4t+2$ for some $t\in\mathbb{Z}$. This imples that $\frac{p-p^2}{2}\equiv1(\operatorname{mod}2)$. Since $\lVert z\rVert\equiv1(\operatorname{mod}4)$, $z_1$ must be odd integer. By equation (\ref{last3}), we have
	\begin{eqnarray*}
		0\equiv2b_1^2+2z_1b_1+2 (\operatorname{mod}4).
	\end{eqnarray*}
	This implies that $2\equiv2(b_1(b_1+z_1))(\operatorname{mod}4)$.
	However, for any integer $b_1,z_1$ with $z_1$ is odd integer, $(b_1(b_1+z_1))\equiv0,2(\operatorname{mod}4)$, i.e., $2\equiv0(\operatorname{mod}4)$, which is a contradiction. Then, there are no integers $b_1,b_2$ corresponding to equation (\ref{last3}). By Theorem \ref{mainthm}, we can conclued that $A(p,z)$ doesn not satisfy the conjecture.
\end{proof}

 For example, when $z=1+2\sqrt{10}$ and $p=3$, we see that the assumptions of Corollary \ref{coroppose} are fulfilled and thus  $
 	A(3,1+2\sqrt{10})$ does not satisfy the Cossu and Zanardo's conjecture.	 Moreover, if we hope that  $$A(3,1+2\sqrt{10})= \begin{pmatrix}
	a_1+a_2\sqrt{10}&b_1+b_2\sqrt{10}\\c_1+c_2\sqrt{10}&1-a_1-a_2\sqrt{10}
\end{pmatrix} \begin{pmatrix}
	d_1+d_2\sqrt{10}&e_1+e_2\sqrt{10}\\f_1+f_2\sqrt{10}&1-d_1-d_2\sqrt{10}
\end{pmatrix},$$
for some $a_1,a_2,b_1,b_2,c_1,c_2,d_1,d_2,f_1,f_2\in\mathbb{Z}$, then by Lemma \ref{lem01},  the system \begin{eqnarray}
	3&=&a_1d_1+10a_2d_2+b_1f_1+10b_2f_2\label{ex3}, \\	0&=&d_1a_2+d_2a_1+b_1f_2+b_2f_1\label{ex0},
\end{eqnarray}
\begin{equation*}
	\begin{split}
			b_1&=\frac{1-a_1-20a_2}{-13},\\
		b_2&=\frac{2-2a_1-a_2}{-13},\\
		c_1&=\frac{a_1-20a_2}{3},\\
		c_2&=\frac{-2a_1+a_2}{3},
	\end{split}\qquad\qquad
	\begin{split}
		e_1&=\frac{d_1+20d_2}{3},\\
		e_2&=\frac{2d_1+d_2}{3},\\
		f_1&=\frac{1-d_1+20d_2}{-13},\\
		f_2&=\frac{-2+2d_1-d_2}{-13},
	\end{split}
\end{equation*}
must have an integral solution.  By using the modulo technique for solving the above system equations, we must have $a_1=39l-7a_2-12$ and $d_1=39n+7d_2+27$, for some $l,n\in\mathbb{Z}$ and for any $a_2,d_2\in\mathbb{Z},$ such that $b_1,b_2,c_1,c_2,e_1,e_2f_1,f_2\in\mathbb{Z}$. Substitute the values $a_1,a_2,d_1,d_2$ into the equation (\ref{ex3}) and (\ref{ex0}), then we obtained a system of Diophantine equations that \begin{eqnarray*}
	3&=&1170nl+210ld_2+819l-210na_2-345n-147a_2-65d_2-242\\
	0&=&30ld_2-12n+21a_2-11d_2+30na_2-8.
\end{eqnarray*}
Surprisingly, we obtain that the above system of equations does not have an integer solution by using an online calculator for Diophantine equations in \cite{web2}. Hence, we can conclude that $A(3,1+2\sqrt{10})$ cannot be written as a product of any two idempotent matrices over $\mathbb{Z}[\sqrt{10}]$.  
In addition, we have applied the same process to many other matrices that do not correspond to the Cossu and Zanardo's conjecture, e.g. $
A(	7,13+2\sqrt{10}),
	A(3,\sqrt{15})$
 and 
$A(	7,\sqrt{77})$.  It turns out that these matrices can not be written as a product of two idempotent matrices, which motivates us to believe that:

\begin{conj} Let $z\in I_p(D)$.
	If $\begin{pmatrix}
		p&z\\ \bar{z}&k
	\end{pmatrix}$ can be written as a product of two idempotent matrices, then $\begin{pmatrix}
		p&z\\ \bar{z}&k
	\end{pmatrix}$ must satisfy the Cossu and Zanardo's conjecture.
\end{conj}

\section*{Acknowledgments}
The authors would like to thank anonymous referee(s) for reviewing this manuscript.  The first author also would like to thank  center of excellence  in non linear analysis and optimization, Naresuan University, for some financial support.  The second author also would like to thank Faculty of Science, Naresuan University,  for financial support on the project number R2566E014.

\end{document}